\documentclass{proc-l}
\usepackage{amssymb,amsmath,amscd}
\usepackage[all]{xy}

\begin{document}
\newtheorem{corollary}{Corollary}
\newtheorem{theorem}[corollary]{Theorem}
\newtheorem*{lemma1}{Lemma}
\newtheorem{lemma}[corollary]{Lemma}

\theoremstyle{definition}
\newtheorem*{definition1}{Definition}
\newtheorem{definition}[corollary]{Definition}

\theoremstyle{remark}
\newtheorem{remark}[corollary]{Remark}

\newcommand{\cC}{\mathcal{C}}
\newcommand{\cD}{\mathcal{D}}
\newcommand{\cE}{\mathcal{E}}
\newcommand{\cF}{\mathcal{F}}
\newcommand{\cG}{\mathcal{G}}
\newcommand{\cH}{\mathcal{H}}
\newcommand{\cM}{\mathcal{M}}
\newcommand{\cN}{\mathcal{N}}
\newcommand{\cO}{\mathcal{O}}
\newcommand{\cS}{\mathcal{S}}
\newcommand{\cV}{\mathcal{V}}
\newcommand{\ccV}{{}^c\mathcal{V}}
\newcommand{\ch}{\mathrm{ch}}
\newcommand{\coker}{\mathrm{coker}}
\newcommand{\cssX}{{}^c\!S^*X}
\newcommand{\ctX}{{}^cTX}
\newcommand{\ctsX}{{}^cT^*X}
\newcommand{\cun}{\cC^{\infty}}
\newcommand{\cuno}{\cC^{\infty}_0}
\newcommand{\cz}{{\mathbb C}}
\newcommand{\dc}{\cD_c}
\newcommand{\Hom}{\mathrm{Hom}}
\newcommand{\iH}{{I_H}}
\newcommand{\iM}{{I_M}}
\newcommand{\ind}{\mathrm{index}}
\newcommand{\nz}{{\mathbb N}}
\newcommand{\pc}{\Psi_{\!c}}
\newcommand{\pcs}{\Psi_{\!c,\sus}}
\newcommand{\ps}{\Psi_\sus}
\newcommand{\px}{\partial_x}
\newcommand{\re}{\mathrm{Re}}
\newcommand{\Res}{\mathrm{Res}}
\newcommand{\rc}{{\rho_c}}
\newcommand{\rz}{{\mathbb R}}
\newcommand{\sD}{\sigma(D)}
\newcommand{\sign}{\mathrm{sign}}
\newcommand{\sssM}{{S_\sus^*(M)}}
\newcommand{\ssX}{S^*X}
\newcommand{\sus}{{\mathrm{sus}}}
\newcommand{\supp}{\mathrm{supp}}
\newcommand{\ta}{\tilde{a}}
\newcommand{\tD}{{\tilde{D}}}
\newcommand{\td}{{\tilde{\delta}}}
\newcommand{\tT}{{\tilde{T}}}
\newcommand{\Td}{\mathrm{Td}}
\newcommand{\Tr}{\mathrm{Tr}}
\newcommand{\tssM}{T_\sus^* M}
\newcommand{\tX}{{\tilde{X}}}
\newcommand{\ucz}{\underline{\cz}}
\newcommand{\urz}{\underline{\rz}}
\newcommand{\zz}{{\mathbb Z}}

\title[On Carvalho's formulation of cobordism invariance]
{On Carvalho's $K$-theoretic formulation of the cobordism invariance
of the index}
\author{Sergiu Moroianu}
\thanks{Partially supported by the
RTN HPRN-CT-2002-00280
``Quantum Spaces -- Noncommutative Geometry'' and the Marie Curie MERG 006375
funded by the European Commission, and by a CERES contract (2004)}
\date{\today}
\address{Institutul de Matematic\u{a} al Academiei Rom\^{a}ne\\
P.O. Box 1-764\\RO-014700
Bucha\-rest, Romania}
\email{moroianu@alum.mit.edu}
\begin{abstract}
We give an analytic proof of the 
fact that the index of an elliptic operator on the boundary of a 
compact manifold vanishes when the principal symbol comes from the 
restriction of a $K$-theory class from the interior. The proof
uses noncommutative residues inside
the calculus of cusp pseudodifferential operators of Melrose.
\end{abstract}
\maketitle

\section*{Introduction}

A classical result of Thom states that the topological signature of 
the boundary of a compact manifold with boundary vanishes. Regarding the 
signature as the index of an elliptic operator, Atiyah and Singer \cite{as}
generalized this vanishing to the so-called twisted signatures.
The cobordism invariance of the index, as this vanishing is known,
was the essential step in their first
proof of the index formula on closed manifolds. 
Conversely, cobordism invariance follows from the index theorem
of \cite{as}. 

On open manifolds a satisfactory index formula 
is not available, and probably not reasonable to expect in 
full generality. Such formulae in various particular cases are 
given e.g., in \cite{aps1}, \cite{melaps} for manifolds with boundary,
in \cite{boris}, \cite{phind} for manifolds with fibered boundary, and in 
\cite{loya}, \cite{in} for manifolds with corners in the sense of Melrose.
To advance in this direction, we believe it is important to understand 
conditions which ensure the vanishing of the index, in particular cobordism
invariance, without using any index formula.

Direct proofs of the cobordism invariance of the index for first-order 
differential operators on closed manifolds
were given e.g., in \cite{braver}, \cite{higson}, \cite{lesch},
\cite{nicol}, and also \cite[Theorem 1]{cii}. We have proposed in 
\cite{cii} an extension of cobordism invariance to manifolds with
corners. The result states that the sum
of the indices on the hyperfaces is null, under suitable hypothesis.

All these results are partial, in that they only 
apply to differential operators of a special type. 
A well-known fact states that the index of ``geometrically defined''
operators is cobordism-invariant; but besides being vague, 
this is also not true (look at the Gau\ss-Bonnet operator).
Only very recently, Carvalho \cite{carvalho,carvart} found a remarkable 
$K$-theoretic statement of cobordism invariance of the index on open manifolds,
using the topological approach of \cite{as2}. 
Here is a reformulation
of the main result of \cite{carvalho} specialized to closed manifolds:
\begin{theorem} \label{th2}
Let $M$ be the boundary of the compact manifold $X$ and $D$ an elliptic 
pseudodifferential operator on $M$. The principal symbol of $D$ defines
a vector bundle over the sphere
bundle inside $T^*M\oplus\urz$. If the class in $K^0(S(T^*M\oplus\urz))$
of this bundle is the restriction of a class from
$K^0(\ssX)$ modulo $K^0(M)$, then $\ind(D)=0$. 
\end{theorem}
The missing details appear in Theorem \ref{kr}.
The aim of this note is to reprove Theorem \ref{th2} with analytic methods.
In order to make the proof likely to generalize to
open manifolds, we have made a point of avoiding 
to use results from $K$-theory, e.g., Bott periodicity 
and the index theorem. Our approach is 
based on Theorem \ref{pt}, a statement about 
the cusp calculus of pseudodifferential 
operators of Melrose on the manifold with boundary $X$, 
in the spirit of \cite{cii}. 

Although they do not appear explicitly in the literature, 
Carvalho's statement (in the closed manifold case) and its present 
variant could be recovered from known results in $K$-theory and $K$-homology. 
We would like to mention here only \cite[Prop. 3]{mepi}, whose arrow-theoretic  
proof could be extended to pseudodifferential operators. For completeness, we 
show in Section \ref{clo} how to retrieve Theorem \ref{th2} also
from the Atiyah-Singer index formula. 

\subsection*{Acknowledgments} I am grateful to Paolo Piazza for useful 
discussions, and to the anonymous referee for valuable 
remarks which greatly improved Section \ref{kre}.

\section{Review of Melrose's cusp algebra} \label{rev}

In this section we recall the facts about the cusp algebra needed in the sequel. 
For a full treatment of the cusp algebra we refer to \cite{meni96c} and 
\cite{in}.

Let $X$ be a compact manifold with boundary $M$, and $x:X\to\rz_+$ a 
boundary-defining function. Choose a product decomposition 
$M\times[0,\epsilon)\hookrightarrow X$. 
A vector field $V$ on $X$ is called \emph{cusp} if
$dx(V)\in x^2\cun(X)$. The space of cusp vector fields forms a 
Lie subalgebra $\ccV(X)\hookrightarrow \cV(X)$ which is a finitely 
generated projective $\cun(X)$-module; indeed, a local basis of $\ccV(X)$ 
is given by $\{x^2\partial_x, \partial_{y_j}\}$ where $y_j$ 
are local coordinates on $M$. Thus there exists a vector bundle $\ctX\to X$ 
such that $\ccV(X)=\cun(X,\ctX)$. Fix a Riemannian metric $g$
on $X\setminus M$ of the form $dx^2/x^4+h^M$ near $x=0$; it extends
to a metric on the fibers of $\ctX$ over $X$ and is 
called a \emph{cusp metric} on $X$ (such $g$ is traditionally called an 
\emph{exact} cusp metric).

The algebra $\dc(X)$ of (scalar) cusp differential operators is defined
as the universal enveloping algebra of $\ccV(X)$ over $\cun(X)$. In 
a product decomposition as above, an operator in $\dc(X)$ of order $m$
takes the form
\begin{equation}\label{cudi}
P=\sum_{j=0}^m P_{m-j}(x) (x^2\px)^j
\end{equation}
where $P_{m-j}(x)$ is a smooth family of differential operators 
of order $m-j$ on $M$.

\subsection{Cusp pseudodifferential operators} The operators in 
$\dc(X)$ can be described alternately (see \cite{meni96c}) 
in terms of their Schwartz kernels.
Namely, there exists a manifold with corners $X^2_c$ obtained by blow-up from 
$X\times X$, and a submanifold $\Delta_c$, such that $\dc(X)$ corresponds
to the space of distributions on $X^2_c$ which are classical conormal to 
$\Delta_c$, supported on $\Delta_c$ and smooth at the boundary face of $X^2_c$
which meets $\Delta_c$. It is then a showcase application of Melrose's
program \cite{meicm} to construct a calculus of 
pseudodifferential operators
$\pc^\lambda(X)$, $\lambda\in\cz$, in which $\dc(X)$ sits as the 
subalgebra of differential operators (the symbols used in the definition 
are classical of order $\lambda)$. No extra difficulty appears
in defining cups operators acting between sections of vector bundles over $X$.
By adjoining the multiplication
operators by $x^z$, $z\in\cz$, we get a pseudodifferential
calculus with two complex indices
\[\pc^{\lambda,z}(X,\cF,\cG):=x^{-z}\pc^\lambda(X,\cF,\cG)\]
such that $\pc^{\lambda,z}(X,\cE,\cF)\subset \pc^{\lambda',z'}(X,\cE,\cF)$ 
if and only if $\lambda'-\lambda\in\nz$ and $z'-z\in\nz$ (since we work 
with classical symbols). Also,
\[\pc^{\lambda,z}(X,\cG,\cH)\circ\pc^{\lambda',z'}(X,\cF,\cG)\subset 
\pc^{\lambda+\lambda',z+z'}(X,\cF,\cH).\]

The fixed cusp metric and a metric on $\cF$ allow one to define the space of cusp 
square-integrable sections $L^2_c(X,\cF)$.
By closure, cusp operators act on a scale of weighted Sobolev spaces 
$x^\alpha H_c^\beta$: \[\pc^{\lambda,z}(X,\cF,\cG)\times 
x^\alpha H_c^\beta(X,\cF)\to x^{\alpha-\Re(z)}H_c^{\beta-\Re(\lambda)}(X,\cG).\]

\subsection{Symbol maps}
There exists a natural surjective \emph{cusp principal symbol} map from 
$\pc^{\lambda}$ onto the space of homogeneous functions on 
$\ctsX\setminus\{0\}$ of homogeneity $\lambda$, which extends the usual 
principal symbol map over the interior of $X$:
\[\sigma:\pc^\lambda(X,\cE,\cF)\to \cun_{[\lambda]}(\ctsX,\cE,\cF).\] 
In the sequel we refer to $\sigma$ as the principal symbol map.
A cusp operator is called \emph{elliptic} if 
its (cusp) principal symbol is invertible on $\ctsX\setminus\{0\}$. 

\begin{definition1}[\cite{meleta}]
Let $\Psi_\sus^\lambda(M,\cE,\cF)$ be the space of classical pseudo-differential
operators $P$ of order $\lambda\in\cz$ from $\cE$ to $\cF$ which are
translation invariant, and such that the convolution kernel 
$\kappa_P(x,y_1,y_2)$ (which is smooth for $x\neq 0$)
decays rapidly as $|x|$ tends to infinity.
\end{definition1}
Under partial Fourier transform in the variable $x$, $\Psi_\sus(M,\cE,\cF)$ 
is identified with the space
of families of operators on $M$ with one real parameter $\xi$, with joint 
symbolic behavior in $\xi$ and in the cotangent variables of $T^*M$.

The second symbol map
is a surjection $I_M:\pc^{\lambda}(X,\cE,\cF)\to 
\Psi_\sus^\lambda(M,\cE,\cF)$, called the \emph{indicial family} map 
\cite{meleta}. If $P$ is given by Eq.\ \eqref{cudi} near $x=0$, then 
\[I_M(P)(\xi)=\sum_{j=0}^m P_{m-j}(x) (i\xi)^j.\]

The principal symbol map and the indicial family are star-morphisms, 
i.e., they are multiplicative and commute with taking adjoints. 
Elliptic cusp operators whose indicial family is invertible for each $\xi\in\rz$
are called fully elliptic. Being fully elliptic is equivalent to being 
Fredholm (see \cite{melaps}).

Let $L^\lambda:=\{(U,\alpha)\in\Psi_\sus^\lambda(M,\cE,\cF)\times
\cun_{[\lambda]}(\ctsX,\cE,\cF); \sigma(U)=\alpha_{|x=0}\}$.
It is proved in \cite{meni96c} that the joint symbol map
\begin{equation}\label{jss}
(\sigma_\lambda,I_M):\pc^{\lambda}(X,\cE,\cF)\to L^\lambda
\end{equation} 
is surjective.

\subsection{Analytic families of cusp operators}\label{afoco}
Let $Q\in\pc^{1,0}(X,\cE)$ be a positive fully elliptic cusp operator of order 
$1$. Then the complex powers $Q^\lambda$ form an analytic family of 
cusp operators of order $\lambda$. 

Let $\cz^2\ni(\lambda,z)\mapsto P(\lambda,z)\in\pc^{\lambda,z}(X,\cE)$ 
be an analytic family in two complex variables. Then $P(\lambda,z)$
is trace-class on $L^2_c(M,\cE)$ for $\Re(\lambda)<-\dim(X), \Re(z)<-1$. 
Moreover, $(\lambda,z)\mapsto\Tr(P(\lambda,z))$ is analytic,
extends to $\cz^2$ meromorphically with at most simple poles 
in each variable at $\lambda\in\nz-\dim(X)$, $z\in\nz-1$, and
\begin{equation}\label{trz}
\Res_{z=-1}\Tr(P(\lambda,z))=\frac{1}{2\pi}\int_\rz \Tr(\iM(x^{-1}P(\lambda,-1)))
d\xi.
\end{equation}
This identity is the content of \cite[Prop.\ 3]{cii}.

\section{Cobordism invariance for cusp operators}\label{cobcus}

This section extends a result from \cite{cii} to 
pseudodifferential operators, in a form which can be applied to $K$-theory. 
We use the same line of proof, with some extra technical difficulties. 
A similar extension from the differential to the 
pseudodifferential case appears
in \cite{kso} when computing the $K$-theory of the algebra $\ps^0(M)$. 
\begin{theorem}\label{pt}
Let $X$ be a compact manifold with boundary $\partial X=M$, and
\[D:\cun(M,\cE^+)\to\cun(M,\cE^-)\]
a classical pseudodifferential operator of order $1$ on $M$. 
Assume that there exist hermitian vector 
bundles $V^+,V^-\to M$, $\cG\to X$ with 
$\cG|_M=\cE^+\oplus\cE^-\oplus V^+\oplus V^-$, and
an elliptic symmetric cusp pseudodifferential 
operator $A\in\pc^{1,0}(X,\cG)$ such that
\begin{equation}\label{noua}
\iM(A)(\xi)=\left[\begin{smallmatrix}
\xi&\tD^*(\xi)&&\\ \tD(\xi)&-\xi &&\\
&&(1+\xi^2+\Delta^+)^{\frac12}&&\\&&&&
-(1+\xi^2+\Delta^-)^{\frac12}
\end{smallmatrix}\right]
\end{equation}
where $\Delta^+,\Delta^-$ are connection Laplacians on $V^+,V^-$,
$\tD\in\ps^1(M,\cE^+,\cE^-)$
and $\tD(0)=D$.
Then $\ind(D)=0$.
\end{theorem}
\begin{proof} 
We first show that we can assume without loss of generality that
$D$ is either injective or surjective. Assuming this, we construct from $A$
a positive cusp operator $Q$ of order $1$. The complex powers of $Q$
are used in defining 
a complex number $N$ as a non-commutative residue. The proof will be finished 
by computing $N$ in two ways; first we get $N=0$, then $N$ is shown to be
essentially $\ind(D)$.

\subsection*{Reduction to the case where $D$ is injective or surjective}
Fix an operator $T\in\Psi^{-\infty}(M,\cE^+,\cE^-)$ such that $D+T$
is either injective or surjective (or both). Choose
$\tT\in\Psi_\sus^{-\infty}(X,\cE^+,\cE^-)$ with $\tT(0)=T$. Choose 
$S\in\pc^{-\infty,0}(X,\cG)$ such that 
\[I_M(S)(\xi)=\begin{bmatrix}
&\tT^*(\xi)&&\\ \tT(\xi)&&&\\&&0&\\&&&0\end{bmatrix}.\]
We can assume that $S$ is symmetric (if not, replace $S$ by $(S+S^*)/2$).
Replace $D$ by $D+T$ and $A$ by $A+S$. Note that $\ind(D)=\ind(D+T)$, since
$T:H^1_c\to L^2_c$ is compact.
The hypothesis of the theorem (in particular (\ref{noua})) still hold for 
$D+T$ instead of $D$ and with $A+S$ instead of $A$. So we can additionally
assume that $D$ is surjective or injective.

\subsection*{Construction of a positive cusp operator $Q$}
For $\xi\in\rz$ we have $\sigma_1(\tD(\xi))=\sigma_1(D)$, so 
$\tD(\xi)$ is elliptic as an operator on $M$ and 
$\ind(\tD(\xi))=\ind(D)$. If $D$ is surjective or injective,
then $0$ does not belong to the spectrum of $DD^*$ (respectively $D^*D$)
so by continuity $\tD(\xi)$ will have the same property for small enough
$|\xi|$. Thus there exists $\epsilon>0$ such that
the kernel and the cokernel of $\tD(\xi)$ have constant 
dimension (hence they vary smoothly)
for all $|\xi|<\epsilon$. Choose a smooth real function 
$\phi$ supported in $[-\epsilon,\epsilon]$ such that $\phi(0)=1$. 
By \cite[Lemma 2]{cii} and the choice of $\phi$, the families 
$\phi(\xi)P_{\ker \tD(\xi)}$ and $\phi(\xi)P_{\ker \tD(\xi)}$
define suspended operators in $\ps^{-\infty}(M)$. 
Let $R\in\pc^{-\infty,0}(X,\cG)$ be such that
\begin{equation}\label{imr}\begin{split}
\iM(R)(\xi)&=\begin{bmatrix}
\phi(\xi)P_{\ker \tD(\xi)}&&&\\
&\phi(\xi)P_{\coker \tD(\xi)}&&\\
&&0&\\&&&0
\end{bmatrix}\\
&\quad\in\Psi_\sus^{-\infty}(M,\cE^+\oplus\cE^-\oplus V^+\oplus V^-).
\end{split}\end{equation}
It follows that $\iM(A^2+R^*R)(\xi)$ is invertible for all $\xi\in\rz$, so 
the cusp operator $A^2+R^*R$ is fully elliptic; this implies that it 
is Fredholm, and moreover its kernel is made of smooth sections 
vanishing rapidly towards $\partial X$. Let $P_{\ker(A^2+R^*R)}$
be the orthogonal projection on the finite-dimensional nullspace
of $A^2+R^*R$. Clearly $A^2+R^*R\geq 0$,
thus $A^2+R^*R+P_{\ker(A^2+R^*R)}$ is strictly positive.
Set 
\[Q:=(A^2+R^*R+P_{\ker(A^2+R^*R)})^{1/2}\]
and let $Q^\lambda$ be the complex powers of $Q$. 
Since $Q^2-A^2\in \pc^{-\infty,0}(X,\cG)$ and $A$ is self-adjoint, we deduce
that for all $\lambda\in\cz$,
\begin{equation}\label{ci}
[A,Q^\lambda]\in\pc^{-\infty,0}(X,\cG).
\end{equation}

\subsection*{A non-commutative residue}
Let $P(\lambda,z)\in\pc^{-\lambda -1, -z-1}(X,\cG)$ be the analytic family 
of cusp operators 
\[P(\lambda,z):=[x^z,A]Q^{-\lambda -1}.\]
From \eqref{trz}, $\Tr(P(\lambda,z))$
is holomorphic in $\{(\lambda,z)\in\cz^2; \Re(\lambda)>\dim(X)-1, \Re(z)>0\}$
and extends meromorphically to $\cz^2$. 
Following the scheme of \cite[Theorem 1]{cii}, our proof of Theorem \ref{kr} 
will consist of computing in two different ways the complex number 
\[N:=\Res_{\lambda=0} \left(\Tr(P(\lambda,z))|_{z=0}\right),\]
i.e., $N$ is the coefficient of $\lambda^{-1}z^0$ in the Laurent expansion 
of $\Tr(P(\lambda,z))$ around $(0,0)$. The idea is to evaluate at $z=0$
\emph{before} and then \emph{after} taking the residue at $\lambda=0$, noting that the final answer
is independent of this order.

\subsection*{Vanishing of $N$}
On one hand, 
\[P(\lambda,z)=x^z[A,Q^{-\lambda -1}]+[A,Q^{-\lambda -1}x^z].\]
The meromorphic function $\Tr[A,Q^{-\lambda -1}x^z]$ is identically zero 
since it vanishes on the open set 
$\{(\lambda,z)\in\cz^2; \Re(\lambda)>\dim(X)-1, \Re(z)>0\}$
by the trace property. By \eqref{ci}, the function 
$\Tr(x^z[A,Q^{-\lambda -1}])$ is regular in $\lambda\in\cz$, so
in particular the meromorphic function
\[z\mapsto\Res_{\lambda=0}\Tr(x^z[A,Q^{-\lambda -1}])\]
vanishes. We conclude that $N=0$. 

\subsection*{Second computation of $N$}
On the other hand, $P(\lambda,0)=0$ so 
\[U(\lambda,z):=z^{-1}P(\lambda, z)\in 
\pc^{-\lambda -1, -z-1}(X,\cG)\]
is an analytic family in $\pc(X,\cG)$. Set
$[\log x,A]:=(z^{-1}[x^z,A])|_{z=0}\in\pc^{0,1}(X,\cG)$. Then 
$U(\lambda, 0)=[\log x,A]Q^{-\lambda-1}$. By multiplicativity of 
the indicial family,
\[\iM(x^{-1}U(\lambda,0))=\iM(x^{-1}[\log x,A])\iM(Q^{-\lambda-1}).\]
By \eqref{noua} and \cite[Lemma 3.4]{in}, 
we see that $\iM(x^{-1}[\log x,A])$ is the $4\times 4$
diagonal matrix
\[\begin{bmatrix}
i&&&\\&-i&&\\&& i\xi(1+\xi^2+\Delta^+)^{-\frac12}&\\&&&
-i\xi(1+\xi^2+\Delta^-)^{-\frac12}
\end{bmatrix}\] 
and  $\iM(Q^{-\lambda-1})=
\iM(A^2+R^*R)^{-\frac{\lambda+1}{2}}$. Also, using \eqref{imr}, we deduce that
$\iM(A^2+R^*R)$ is the $4\times 4$ diagonal matrix with entries 
\begin{align*}
a_{11}&=\xi^2+\tD(\xi)^*\tD(\xi)+\phi^2(\xi)P_{\ker \tD(\xi)}&
a_{33}&=1+\xi^2+\Delta^+\\
a_{22}&=\xi^2+\tD(\xi)\tD(\xi)^*+\phi^2(\xi)P_{\coker \tD(\xi)}&
a_{44}&=1+\xi^2+\Delta^-.
\end{align*}
By \eqref{trz},
\begin{equation*}\begin{split}
\Tr(P(\lambda, z))|_{z=0}
&=\frac{1}{2\pi}\int_\rz \Tr(\iM(x^{-1}(U(\lambda,0)))d\xi\\
&=\frac{i}{2\pi}\int_\rz 
\left(\Tr(\xi^2+\tD(\xi)^*\tD(\xi)
+\phi^2(\xi)P_{\ker \tD(\xi)})^{-\frac{\lambda+1}{2}}\right.\\
&\quad-\Tr(\xi^2+\tD(\xi)\tD(\xi)^*
+\phi^2(\xi)P_{\coker \tD(\xi)})^{-\frac{\lambda+1}{2}}\\
&\quad\left.+\xi\Tr(1+\xi^2+\Delta^+)^{-\frac\lambda2-1}
-\xi\Tr(1+\xi^2+\Delta^-)^{-\frac\lambda2-1}\right)
d\xi
\end{split}\end{equation*}
The third and fourth terms in this sum are odd in $\xi$ 
so their integral vanishes. For each fixed $\xi$ we compute
the trace of the first two terms by using orthonormal basis
of $L^2_c(M,\cE^+)$, $L^2_c(M,\cE^-)$ given by eigensections of
$\tD(\xi)^*\tD(\xi)$, respectively $\tD(\xi)\tD(\xi)^*$. The non-zero parts
of the spectrum of $\tD(\xi)^*\tD(\xi)$ and $\tD(\xi)\tD(\xi)^*$ coincide,
so what is left is 
\[\int_\rz \ind(\tD(\xi))
(\xi^2+\phi^2(\xi))^{-\frac{\lambda+1}{2}}d\xi.\]
The subtle point here is that the kernel and cokernel of $\tD(\xi)$ 
may have jumps when $|\xi|>\epsilon$, but our formula involves only 
the index, which is homotopy invariant and equals $\ind(D)$ for all $\xi\in\rz$.
Thus the index comes out of the integral; the residue
\[\Res_{\lambda=0}\int_\rz (\xi^2+\phi^2(\xi))^{-\frac{\lambda+1}{2}}d\xi\]
is independent of the compactly supported function $\phi$ and equals $2$, so
\[0=N=\Res_{\lambda=0}\Tr(P(\lambda, z))|_{z=0}=\frac{i}{\pi}\ind(D).\qedhere\]
\end{proof}

\section{The $K$-theoretic characterization of cobordism invariance}\label{kre}

We interpret now Theorem \ref{pt} in topological terms. 
Let 
\[p:\sssM\to M\] be the sphere bundle inside $\tssM:=T^*M\oplus\underline{\rz}$.
We also denote by $p$ the bundle projections $T^*M\to M$, $T^*X\to X$, 
$S^*X\to X$. The total space of $\sssM$ is the boundary of $\cssX$. By fixing 
a product decomposition of $X$ near $M$, we get non-canonical vector bundle 
isomorphisms making the diagram
\begin{equation*}
\xymatrix{
\ctsX\ar[d]^\cong\ar[r]^r&T^*_\sus M\ar[d]^\cong\\
T^*X\ar[r]^r&T^*X|_M}
\end{equation*}
commutative, so we can replace $\cssX$ with the more familiar space 
$\ssX$ in all the topological considerations of this section.

The interior unit normal vector inclusion $\imath:M\to\sssM$ and the 
bundle projection $p:\sssM\to M$ induce a splitting
\[K^0(\sssM)=\ker(\imath)\oplus p^*(K^0(M)).\]
Let 
\[r:K^0(\ssX)\to K^0(\sssM)\]
be the map of restriction to the boundary, and
\[d:K^0(T^*M)\to K^0(\sssM)/ p^*(K^0(M))\]
the isomorphism defined as follows: if 
$(\cE^+,\cE^-,\sigma)$ is a triple defining a class in $K^0(T^*M)$ with 
$\sigma:\cE^+\to\cE^-$ an isomorphism outside the open unit ball,
then 
\[
d(\cE^+,\cE^-,\sigma)=\begin{cases}\cE^+ &\text{on
$\sssM\cap\{\xi\geq 0\}$}\\
\cE^- &\text{on
$\sssM\cap\{\xi\leq 0\}$}
\end{cases}\]
with the two bundles identified via $\sigma$ over 
$\sssM\cap\{\xi=0\}= S^*M$. 
We can now reformulate Theorem \ref{th2} as follows:
\begin{theorem}\label{kr}
Let $X$ be a compact manifold with closed boundary $M$, $\cE^\pm\to M$ 
hermitian vector bundles, and $D\in\Psi(M,\cE^+,\cE^-)$ 
an elliptic pseudodifferential operator with symbol class
\[[\sD]:=(p^*\cE^+,p^*\cE^-,\sigma(D))\in K^0(T^* M)\] 
Assume that $d[\sD]\in p^*(K^0(M)) + r(K^0(\ssX)).$
Then $\ind(D)=0.$
\end{theorem}
\begin{proof} 
The idea is to construct an operator $A$ as in Theorem \ref{pt}.
We must first construct the vector bundles $V^\pm$, and then 
extend the principal symbol of \eqref{noua} to an elliptic symbol in 
the interior of $X$. Note that none of the bundles $\cE^\pm, V^\pm$ 
has any reason to extend to $X$.

We can assume that $D$ is of order $1$. Extend $\sD_{|S^*M}$ arbitrarily 
to a homomorphism $\sigma:p^*\cE^+\to p^*\cE^-$ (not necessarily invertible) 
over $\sssM$. Let $\cF^\pm\to\sssM$ be the vector bundles 
defined as the span of the eigenvectors of positive, resp.\ negative 
eigenvalues of the symmetric automorphism of $p^*(\cE^+\oplus\cE^-)$
\[a:=\begin{bmatrix}\xi&\sigma^*\\\sigma&-\xi\end{bmatrix}:
p^*(\cE^+\oplus\cE^-)\to p^*(\cE^+\oplus\cE^-).\]

\begin{lemma1}
The $K$-theory class of the vector bundle $\cF^+$ is $d[\sD]$.
\end{lemma1}
\begin{proof} $\cF^+$ is the image of the projector $\frac{1+a(a^2)^{-\frac12}}{2}$
inside $p^*(\cE^+\oplus\cE^-)$, or equivalently the image of the endomorphism
$(a^2)^{\frac12}+a$:
\begin{equation*}\begin{split}
\cF^+&=\{((\xi+(\xi^2+\sigma^*\sigma)^{\frac12})v,\sigma v) ; v\in\cE^+\}\\
&\quad+\{(\sigma^* w,(-\xi+(\xi^2+\sigma\sigma^*)^{\frac12})w) ; w\in\cE^-\}.
\end{split}\end{equation*}
Now $\xi+(\xi^2+\sigma^*\sigma)^{\frac12}$ is invertible when
$\xi\geq 0$, and $-\xi+(\xi^2+\sigma\sigma^*)^{\frac12}$ is invertible when
$\xi\leq 0$. Thus the projection from $\cF^+$ on $p^*\cE^+$, respectively on 
$p^*\cE^-$, are isomorphisms for $\xi\geq 0$, respectively for $\xi\leq 0$.
Over $\{\xi=0\}$ these isomorphisms differ by $\sigma(\sigma^*\sigma)^{-\frac12}$,
which is homotopic to $\sigma$ by varying the exponent from $-\frac12$ to $0$.
\end{proof}

The hypothesis says therefore that 
\begin{equation}\label{ep}
[\cF^+]\in p^*(K^0(M)) + r(K^0(\ssX)).
\end{equation}

\begin{lemma1}
There exist vector bundles
$G^\pm\to \ssX$, $V^\pm\to M$ such that
\begin{equation}\label{fkg}
\cF^\pm\oplus p^*V^\pm=G^\pm_{|\sssM}
\end{equation}
and moreover there exists $N\in\nz$ with
\begin{align}\label{gtri}\cE^+\oplus \cE^-\oplus V^+\oplus V^-\cong
\ucz^N,&& G^+\oplus G^-\cong\ucz^N.\end{align} 
\end{lemma1}
\begin{proof} From eq.\ \eqref{ep}, there exist $V^+_0\to M$, $G^+_0\to X$ and
$k\in \nz$ with 
$\cF^+\oplus \ucz^k=p^*V^+_0\oplus {G^+_0}_{|\sssM}\oplus\ucz^k$. Let $V^+_1$ 
be a complement of $V^+_0$ inside some $\ucz^h$. Then $V^+:=\ucz^k\oplus V^+_1$
and $G^+:=\ucz^{h+k}\oplus G^+_0$ satisfy \eqref{fkg}. This implies
\begin{equation*}\begin{split}
[\cF^-]&=[p^*(\cE^+\oplus \cE^-)]-[\cF^+]
=p^*[\cE^+\oplus \cE^-\oplus V^+]-r[G^+].
\end{split}\end{equation*}
Let $G^-_0$ be a complement (inside $\ucz^{N_0}$) of $G^+$, and 
$V^-_0$ a complement of $\cE^+\oplus \cE^-\oplus V^+$ inside $\ucz^{N_1}$.
Then \[[\cF^-]+p^*[V^-_0]+\ucz^{N_0}=\ucz^{N_1}+r[G^-_0]\]
which amounts to saying that there exists $N_2\in\nz$ with
\[\cF^-\oplus p^*V^-_0\oplus \ucz^{N_0+N_2}\cong\ucz^{N_1+N_2}
\oplus G^-_0|_\sssM.\]
Thus $V^-:=V^-_0\oplus \ucz^{N_0+N_2}$ and 
$G^-:=\ucz^{N_1+N_2}\oplus G^-_0$ satisfy \eqref{fkg}. 
From the construction of $V^-$ and $G^-$,
Eq.\ \eqref{gtri} holds for $N:=N_0+N_1+N_2$. 
\end{proof}
Let $\cG:=\ucz^N\to X$ be the trivial bundle. 
From \eqref{gtri}, $\cG_{|M}\cong\cE^+\oplus \cE^-\oplus V^+\oplus V^-$
(as bundles over $M$) and $p^*\cG\cong G^+\oplus G^-$ (as bundles over $S^*X$).
Define
$\ta:p^*\cG\to p^*\cG$ to be the automorphism of $p^*\cG$ over $\ssX$
that equals $\pm 1$ on $G^\pm$. From the definition
of $\cF^\pm$ and eq.\ \eqref{fkg} it follows that $\ta_{|\sssM}$ and 
the automorphism $\left[\begin{smallmatrix}a&&\\&1&\\&&-1\end{smallmatrix}
\right]$ (written in the decomposition
$p^*\cG_{|\sssM}=p^*(\cE^+\oplus\cE^-)\oplus p^*V^+\oplus p^*V^-$) have the same spaces of 
eigenvectors of positive, respectively negative eigenvalues. Thus we can deform
$\ta$ inside self-adjoint automorphisms to an automorphism $\alpha$
with
\begin{equation}\label{tas}
\alpha_{|\sssM}=\left[\begin{smallmatrix}a&&\\&1&\\&&-1\end{smallmatrix}\right].
\end{equation}
We extend $\alpha$ to $T^*X\setminus 0$ with homogeneity $1$.

As noted at the beginning of this section, we replace $\ssX$ by $\cssX$.
By (\ref{tas}) and the definition of $a$, $\alpha|_\sssM$ coincides with 
the principal symbol of the right-hand side of (\ref{noua}). Therefore 
using \eqref{jss},
there exists an elliptic cusp operator $A\in\pc^1(X,\cG)$ with 
$\sigma_1(A)=\alpha$ and with 
indicial family given by the symmetric suspended operator (\ref{noua}). 
By replacing $A$ with $(A+A^*)/2$ we can assume 
$A$ to be symmetric. The hypothesis of Theorem \ref{pt} is fulfilled, 
so we conclude that $\ind(D)=0$.
\end{proof}

\section{Variants of Theorem \ref{kr}}\label{clo}
\subsection{Carvalho's theorem}
Carvalho \cite{carvalho} obtained a slightly different statement of 
cobordism invariance (her result holds for non-compact manifolds as well). 
Namely, in the context of Theorem \ref{kr} she proved that 
$\ind(D)=0$ provided that $[\sD]$ lies in the image of 
the composite map
\[K^1(T^*X)\stackrel{r}{\to} K^1(T^*M\oplus \rz) 
\stackrel{\beta^{-1}}{\to} K^0(T^*M)\]
defined by restriction and by Bott periodicity. 
Consider the relative pairs 
\begin{align*}
S^*X\hookrightarrow & B^*X,& \sssM \hookrightarrow & B^*_\sus M
\end{align*} 
the inclusion map between them and the 
induced boundary maps in the long exact sequences in $K$-theory. 
We claim that we get a commutative diagram

\begin{align*}
\xymatrix{
K^0(\ssX)\ar[d]^r\ar @{->>} [r]&K^1(T^*X)\ar[d]^r\\
K^0(\sssM)\ar[d]^q\ar @{->>} [r]&K^1(\tssM)\\
K^0(\sssM)/p^*K^0(M)\ar[r]^-{d^{-1}}&K^0(T^*M)\ar[u]_\beta}\label{dia}
\end{align*}
Indeed, the upper square commutes by naturality and the lower one
by checking the definitions. Moreover, the existence\footnote{The obstruction
lives in $H^n(X)$ which is $0$ when $X$ has nonempty boundary; 
I am grateful to Gustavo Granja for this argument.} 
of nonzero sections 
in $T^*X\to X$ and $\tssM\to M$ shows that the rows are surjective. Also
$\beta$, $d$ are isomorphisms, so $d[\sD]$ lies in the image of $q\circ r$ 
if and only if $[\sD]$ lies in the image of $\beta^{-1}\circ r$.
Thus Theorem \ref{kr}
is equivalent to Carvalho's statement applied to closed manifolds.
Our formulation is marginally simpler because it does not involve 
the Bott isomorphism.

\subsection{An indirect proof of Theorem \ref{kr}} 

As mentioned in the introduction, Theorem \ref{kr} follows from the 
Atiyah-Singer formula:
\[\ind(D)=\langle M, \Td(TM)\cup p_*\ch([\sD])\rangle,\]
where $p_*$ denotes integration along the fibers of $p:T^*M\to M$, taking 
values in the cohomology of $M$ twisted with the orientation bundle.
Indeed, the normal bundle to $M$ in $X$ is trivial so $\Td(TM)=\Td(TX)_{|M}$.
We can embed $T^*M$ into $S^*X_{|M}$ via the central 
projection from the interior pole of each sphere; the pull-back through this 
map of $d[\sD])$ coincides with $[\sD]$ modulo $p^*K^0(M)$, in particular
the push-forward on $M$ of $\ch(d[\sD])$ and of $\ch([\sD])$ are equal.
So the hypothesis that $d[\sD]$ is the restriction of a class on $\ssX$ modulo
$p^*K^0(M)$ implies, by the functoriality of the Chern character, that
$p_*\ch([\sD])\in H^*(M,\cO)$ is the restriction of a (twisted)
cohomology class from $X$. Finally Stokes formula shows that $\ind(D)=0$. 

\bibliographystyle{plain}

\begin{thebibliography}{1}

\bibitem{aps1}
M.F.~Atiyah, V.K.~Patodi and I.M.~Singer,  
{\sl Spectral asymmetry and Riemannian geometry. I,} 
Math. Proc. Cambridge Philos. Soc. {\bf 77} (1975), 43--69.

\bibitem{as}
M.~F.~Atiyah and I.~M.~Singer,
{\sl The index of elliptic operators on compact manifolds, } 
Bull. Amer. Math. Soc. {\bf 69} (1963), 422--433.

\bibitem{as2}
M.~F.~Atiyah and I.~M.~Singer,
{\sl The index of elliptic operators. I,} 
Ann. of Math. {\bf 87} (1968), 484--530.

\bibitem{braver}
M.~Braverman,
{\sl New proof of the cobordism invariance of the index, }
Proc. Amer. Math. Soc. {\bf 130} (2002), 1095--1101.

\bibitem{carvalho}
C.~Carvalho,
{\sl Pseudodifferential operators and applications to Index Theory on 
noncompact manifolds, }
PhD thesis, Trinity College, University of Oxford (2003).

\bibitem{carvart}
C.~Carvalho,
{\sl A $K$-theory proof of the cobordism invariance of the index, }
$K$-Theory {\bf 36} (2005), 1--31.

\bibitem{higson}
N.~Higson,
{\sl A note on the cobordism invariance of the index, }
Topology {\bf 30} (1991), 439--443.

\bibitem{in}
R.~Lauter and S.~Moroianu,
{\sl The index of cusp operators on manifolds with corners, }
Ann. Global Anal. Geom. { \bf 21} (2002), 31--49. 

\bibitem{phind}
R.~Lauter et S.~Moroianu,
{\sl An index formula on manifolds with fibered cusp ends, }
J. Geom. Analysis. {\bf 15} (2005), 261--283.

\bibitem{lesch}
M.~Lesch,
{\sl Deficiency indices for symmetric Dirac operators on manifolds 
with conic singularities, }
Topology {\bf 32} (1993), 611--623.

\bibitem{loya}
P.~Loya,
{\sl Tempered operators and the heat kernel and complex powers of 
elliptic pseudodifferential operators, }
Comm. Partial Differential Equations {\bf 26} (2001), 1253--1321.

\bibitem{melaps}
R.~B.~Melrose,
{\sl The Atiyah-Patodi-Singer index theorem, }
Research Notes in Mathematics {\bf 4}, 
A. K. Peters, Wellesley, MA (1993).

\bibitem{meicm}
R.~B.~Melrose,
{\sl Pseudodifferential operators, corners and singular limits, }
in Proc. Int. Congress of Mathematicians,
  Kyoto, Springer-Verlag Berlin - Heidelberg - New York (1990), 217--234.

\bibitem{meleta}
R.~B.~Melrose,
{\sl The eta invariant and families of pseudodifferential operators, }
Math. Res. Letters {\bf 2} (1995), 541--561.

\bibitem{meni96c}
R.~B.~Melrose and V.~Nistor,
{\sl Homology of pseudodifferential operators {I}. {M}anifolds with
boundary, }
preprint funct-an/9606005.

\bibitem{mepi}
R.~B.~Melrose and P.~Piazza,
{\sl Families of Dirac operators, boundaries and the $b$-calculus, }
J. Differ. Geom. {\bf 46} (1997), 99--180.

\bibitem{kso}
S.~Moroianu,
{\sl $K$-Theory of suspended pseudodifferential operators, } 
$K$-Theory {\bf 28} (2003), 167--181.

\bibitem{cii}
S.~Moroianu,
{\sl Cusp geometry and the cobordism invariance of the index, }
Adv. Math. {\bf 194} (2005), 504--519.

\bibitem{nicol}
L.~I.~Nicolaescu,
{\sl On the cobordism invariance of the index of Dirac operators, }
Proc. Amer. Math. Soc. {\bf 125} (1997), 2797--2801.

\bibitem{boris}
B.~Vaillant,
{\sl Index- and spectral theory for manifolds with generalized
  fibered cusps, }
Dissertation, Bonner Mathematische Schriften {\bf 344} (2001), 
Rheinische Friedrich-Wilhelms-Universit\"at Bonn.

\end{thebibliography}

\end{document}